\UseRawInputEncoding
\documentclass[12pt]{article}
\usepackage{amsfonts}
\usepackage{t1enc}
\usepackage{amsmath}
\usepackage{amssymb}
\usepackage{epsfig}
\usepackage[utf8]{inputenc}							
\usepackage[french,english]{babel}
\usepackage{mathrsfs}
\usepackage{hyperref}
\usepackage{amsthm}
\usepackage{graphics}
\usepackage{layout}
\usepackage{latexsym}
\usepackage[dvipsnames]{xcolor}
\usepackage{indentfirst} 						
\usepackage{bigints}

\usepackage[a4paper]{geometry}
\def\mylabel#1{\label{#1}}

\newtheorem{theorem}{Theorem}[section]
\newtheorem{lemma}[theorem]{Lemma}
\newtheorem{corollary}[theorem]{Corollary}
\newtheorem{proposition}[theorem]{Proposition}
\newtheorem{example}[theorem]{Example}
\newtheorem{examples}[theorem]{Examples}
\newtheorem{remark}[theorem]{Remark}

\newtheorem{definition}[theorem]{Definition}

\newtheorem{assumption}[theorem]{Assumption}

\def\bit{\begin{itemize}}

\def\eit{\end{itemize}}

\reversemarginpar   

\def\bc{\begin{center}}

\def\ec{\end{center}}

\def\bthm{\begin{theorem}}

\def\ethm{\end{theorem}}

\def\bcor{\begin{corollary}}

\def\ecor{\end{corollary}}

\def\bprop{\begin{proposition}}

\def\eprop{\end{proposition}}

\def\blem{\begin{lemma}}

\def\elem{\end{lemma}}

\def\bex{\begin{example}}

\def\eex{\end{example}}

\def\bexs{\begin{examples}}

\def\eexs{\end{examples}}

\def\bexo{\begin{exercice} \rm }

\def\eexo{\end{exercice} }

\def\brem{\begin{remark}}

\def\erem{\end{remark}}

\def\prf{{\bf Proof .}}

\def\bdes{\begin{description}}

\def\edes{\end{description}}

\def\ita{\item[(a)]}

\def\itb{\item[(b)]}

\def\itc{\item[(c)]}

\def\beq{\begin{equation}}

\def\eeq{\end{equation}}

\def\ben{\begin{enumerate}}

\def\een{\end{enumerate}}

\def\beqar{\begin{eqnarray}}

\def\eeqar{\end{eqnarray}}

\def\beqarr{\begin{eqnarray*}}

\def\eeqarr{\end{eqnarray*}}

\def\prf{{\bf Proof. }\hspace{.1in}}

\def\PP{{\mathcal P}}
\def\E{{\mathbb E}}

\def\P{{\mathbb P}}

\def\Ind{{\mathbf 1}}
\def\RR{{\mathbb R}}  

\def\NN{{\mathbb N}}

\begin{document}
\title{Persistence in randomly switched Lotka-Volterra food chains}
\author{Antoine Bourquin\footnote{Email address : antoine.bourquin@unine.ch}
\bigskip
\\ \normalsize  Institut de Math\'ematiques, Universit\'e de Neuch\^atel, Switzerland}


\maketitle

\begin{abstract}
We consider a dynamical system obtained by the random switching between $N$ Lotka-Volterra food chains. Our key assumption will be that at least two vector fields only differ on the resources allocated to the growth rate of the first species. We will show that the existence of a positive equilibrium of the average vector field is equivalent to the persistence of all species. Under this condition, the semi-group converges exponentially quickly to a unique invariant probability measure on the positive orthant. If this condition fails to hold, we have two possibilities. The first possibility is the extinction case, in which a group of species becomes extinct exponentially quickly
while the distribution of the remaining species converges weakly to another invariant probability measure. The second possibility is the critical case, in which there is a weaker form of persistence of some species, whilst some of the remaining become extinct exponentially quickly.
 We will also analyse the sensitivity of this model to the parameters.

\end{abstract}

\paragraph{Keywords} Lotka-Volterra food chains; random switching; piecewise deterministic \\Markov processes; prey-predator; Hörmander condition; stochastic persistence


\tableofcontents



\section{Introduction} \mylabel{sec:intro}


We consider a population model obtained by the random switching between finitely many vector fields. Each environment models a food chain that has the property that species $n$ is the apex predator and hunts species $n-1$, which in turn is the predator of species $n-2$ and so on up to species $1$ which is at the bottom of the food chain.

More precisely, let $E := \{1,\ldots,N\}$ denote the set of possible environments. For each $j\in E$, the dynamic in environment $j$ is given by a Lotka-Volterra food chain differential equation on the state space $\RR_+^n := \left\{x \in \RR^n \mid x_i \geq 0, \; i = 1,\ldots,n \right\}$ defined as
\begin{align} \label{model_by_evironment}
dx(t) = G^j(x(t))dt.
\end{align}
Here $G^j : \RR^n \to \RR^n$ is the vector field defined by 
\begin{align}\label{functions_G^j}
G^j_i(x) := x_i F_i^j(x), \qquad i=1,\ldots,n,
\end{align}
with
\begin{align} \label{fucntions_F_i^j}
F_i^j(x) := \left\{
     \begin{array}{ll}
      a_{10}^j -a_{11} x_1 - a_{12} x_2 								& i=1,\\
     -a_{i0}^j + a_{i,i-1} x_{i-1} - a_{ii} \, x_i - a_{i,i+1} x_{i+1} 	& i= 2, \ldots,n-1, \\
      -a_{n0}^j + a_{n,n-1}  x_{n-1} -a_{nn} \,x_n 					& i=n,
     \end{array}
     \right. 
\end{align}
where $a_{i0}^j,a_{11},a_{ik} >0$ for $i \neq k$, $k\neq0$ and $a_{ii} \geq 0$ for $i=1, \ldots, n$. 

The quantity $a_{ii}$ represents the intra-specific competition rate of species $i$, $a_{i0}^j$ its death rate (or growth rate for species $1$), $a_{i,i-1}$ the rate at which it hunts species $i-1$, and $a_{i,i+1}$ its rate of being hunted by species $i+1$.

Another interpretation of the rate $a_{i0}^j$ is a measure of the habitability of the environment $j$ for species $i$. Indeed, if $a_{i0}^j$ is small (or large for $i=1$) then environment $j$ is good for species $i$ whilst if $a_{i0}^j$ is large (or small for $i=1$) then environment $j$ is inhospitable for species $i$. An example of switching which would only affect these rates would be cleaning up a polluted area, thereby making the environment more liveable.

The switching model we consider in this paper is defined on $\RR_+^n\times E$ by
\begin{align}\label{model}
\left\{ \begin{array}{ll}
dX(t) = G^{J(t)}(X(t)) dt, \\
\P\left(J(t+s) = j \mid \mathcal{F}_t, \, J(t) = i\right) = b_{ij}s + o(s) & \text{for } i\neq j,
\end{array}
\right.
\end{align}
where $J(t)$ is a continuous jump process on $E$, $\mathcal{F}_t = \sigma \{ (X(s),J(s)) \mid s \leq t\}$ is the natural filtration and $b_{ij}> 0$ are the jump rates. The quantity $X_i(t)$ represents the density of species $i$ at time $t$ and $J(t)$ the switching between environments. We let $(Z^{(x,j)}(t))_{t\geq 0}$ denote the solution to \eqref{model} with initial condition $Z^{(x,j)}(0) = (x,j)$.

Let $\varphi^j = \{\varphi^j_t\}$ denote the flow induced by \eqref{model_by_evironment}. We say that $M \subset \RR_+^n$ is {\it positively invariant under $\varphi^j$} if 
$$\varphi_t^j (M) \subset M \qquad \text{for all } t \geq 0.$$
We also say that $M \times E$ is {\it positively invariant for \eqref{model}} if $M$ is positively invariant under $\varphi^j$ for all $j \in E$.

\begin{proposition}\label{compact_B}
There exists a compact set $B \subset \RR_+^n$ such that $B \times E$ is positively invariant for \eqref{model}. Moreover for every $(x,j)\in \RR_+^n\times E$, with probability $1$ there exists a time $t\geq 0$ such that $Z^{(x,j)}(t) \in B \times E$.
\end{proposition}

The proof is postponed to Section \ref{proof_compact_B}.  
This proposition allows us to restrict the state space of the process to $B \times E$ instead of $\RR_+^n \times E$.

This type of process is a particular example of piecewise deterministic Markov process (PDMP) in the sense of Davis \cite{davis}. We refer the reader to \cite{B2018,BMZIHP,BS17} for some general results on PDMP and \cite{BL16,HeningStrickler} for other examples of randomly switching Lotka-Volterra vector fields.

Lotka-Volterra food chains were first considered by T. Gard and T. Hallam \cite{gard1979persistence} in 1979 in the deterministic case with no intra-specific competition. They gave a criterion based uniquely on the coefficients
for persistence and extinction. Recently, A. Hening and D. Nguyen \cite{HeningDang18c,HeningDang18b} extended the model to the stochastic differential equations setting and also considered the case with intra-specific competition. They provided conditions for persistence and extinction, mainly under the assumption that the noise is non-degenerate.
M. Benaim, A. Bourquin and D. Nguyen \cite{LV_degenerate} extended these results to the degenerate situation where the noise affects at least the top or the bottom species.  They also provided some criteria for polynomial rates of convergence.
This paper completes these studies.

We focus here on the situation where there are two environments for which all the rates except the first one ($a_{10}$) are equal. In other words, the switching between these environments only affects the growth rate of the first species.
For instance, this may be the addition of some pollutant or some nutrient only affecting micro-organisms which are at the bottom of the food chain - see \cite{faithfull2011bottom,moriarty1997role} for example.
This situation is in fact very natural. Indeed, natural changes in the wild are generally not abrupt and the first species is usually the most sensitive to small perturbations.  


More precisely, throughout the paper we impose the following assumption.

\begin{assumption}\label{Assumption_bio}
There exist $\beta_1,\beta_2 \in E$ such that 
$$\left\{ \begin{array}{ll}
a_{10}^{\beta_1} \neq a_{10}^{\beta_2}, \\
a_{i0}^{\beta_1} = a_{i0}^{\beta_2} & \text{ for } i=2,\ldots,n.
\end{array} \right. $$
\end{assumption}
%
%

By standard arguments on finite Markov chains, there exists a unique invariant probability measure of $J$ on $E$ defined by
\begin{align}\label{Proba_inv_on_E}
\nu = (\nu_1, \ldots, \nu_N)
\end{align}
where $\nu_j>0$ and $\sum \nu_j =1$. Let
\begin{align}\label{G_nu}
G^\nu := \sum_{j=1}^N \nu_j G^j
\end{align} 
denote the {\it average vector field}.

Roughly speaking, our main results can be summarized as the equivalence of the following conditions:
\begin{enumerate}
\item $G^\nu$ has a positive equilibrium $q^*$ (i.e. $G^\nu(q^*)=0$ and $q^*_i>0$ for $i=1,\ldots,n$).
\item The trajectories induced by $\dot{x} = G^\nu(x)$ converge to $q^*$.
\item The PDMP defined by \eqref{model} is stochastically persistent.
\end{enumerate}

\paragraph*{Outline of contents} Section \ref{main_result} introduces some notation and states the three main results. The first result (Theorem \ref{theorem_stoch_pers}) gives conditions ensuring persistence of all species, the second (Theorem \ref{theorem_extinction}) is devoted to the extinction case  and the third one (Theorem \ref{theorem_critical_case}) deals with the critical case.
In Section \ref{Sect_Math_tolls_stoch_pers} we introduce some general mathematical tools and etablish some properties of the so-called {\it invasion rate} (defined by \eqref{Invasion_rate_equ}). In Section \ref{Section proof} we prove the main theorems and Proposition \ref{compact_B}. Finally, in Section \ref{Section_Sensibility} we investigate the sensitivity of the model subject to small disturbances.


\section{Notation and results} \label{main_result}

\subsection{Notation}

Throughout the paper we let 
$$B_+ := \left\{ x \in B \;\, \middle| \,\; \prod_{i=1}^n x_i >0 \right\} \quad\text{ and }\quad \partial B := \left\{ x \in B \;\, \middle| \,\; \prod_{i=1}^n x_i = 0 \right\}$$
denote respectively the interior and the boundary of $B$. For $1\leq k<n$, we also let
$$B^k := \left\{ x \in B \mid x_1,\ldots,x_k \geq 0 \text{ and } x_{k+1},\ldots,x_n =0 \right\},$$
$$B_+^k := \left\{ x \in B \mid x_1,\ldots,x_k >0 \text{ and } x_{k+1},\ldots,x_n =0 \right\},$$
$$\partial B^k := \left\{ x \in B \;\, \middle| \,\; \prod_{i=1}^k x_i = 0 \right\}.$$

For every initial condition $(x,j) \in B \times E$, we let $(P_t)_{t \ge 0}$ denote the {\it transition kernel} of the process $Z^{(x,j)}(t)$ (the solution to \eqref{model}) defined by $$P_t f(x,j) := \E [f(Z^{(x,j)}(t))]$$ for every measurable bounded function $f : B\times E \to \RR$ and by $$P_t((x,j),A) := P_t \Ind_A (x,j)$$ for all Borel set $A \subset B \times E$. We also let $\P_{(x,j)}$ and $\E_{(x,j)}$ denote the probability measure (respectively the expected value) conditioned on $Z^{(x,j)}(0) = (x,j)$.

A Borel set $A  \subset B \times E$ is called {\it invariant} if $$P_t \Ind_A = \Ind_A, \qquad \text{ for all } t \geq 0.$$

Given an invariant set $A \subset B\times E$ (typically $A= B_+ \times E$ or $A = \partial B \times E$), we let $\PP_{inv}(A)$ (resp. $\PP_{erg}(A)$) denote the {\it set of invariant} (resp. {\it ergodic}) {\it probability measures} of the process that are supported by $A$, i.e. $\mu(A)=1$ for $\mu \in \PP_{inv}(A)$ (resp. for $\mu \in \PP_{erg}(A)$).  Recall that an invariant probability measure $\mu$ is {\it ergodic} if  for every invariant Borel set $B$, $\mu(B)\in\{0,1\}$.

We also let $\left(\Pi_t^{(x,j)}\right)_{t \in \RR_+}$ denote the {\it set of empirical occupation measures} of the process $(Z^{(x,j)}(t))$ with initial condition $(x,j) \in B \times E$, which are defined by
$$\Pi_t^{(x,j)} (\cdot) := \frac{1}{t} \int_0^t \Ind_{\{ Z^{(x,j)}(s) \in \,\cdot \,\}} ds.$$

For a measurable function $f:B \times E \to \RR$ and a measure $\mu$, we write $\mu f := \int f(x) \mu(dx)$ whenever it makes sense.

We recall that the {\it total variation distance} between two probability measures $\mu$ and $\nu $ on $B \times E$ is defined by
$$\|\mu - \nu \|_{TV} := \sup\{|\mu f - \nu f| \mid f:B \times E \to \RR \text{ measurable bounded}, \|f\|_\infty <1\}.$$
%

We say that a sequence of probability measures $(\mu_n)$ on $B \times E$ {\it converges weakly} to a probability measure $\mu$ if for all continuous bounded functions $f : B \times E \to \RR$, $\lim\limits_{n \to \infty}\mu_n f =\mu f$.

\subsection{Main results}

We will start by defining some general notions of persistence and extinction that will be used in this paper:

\begin{enumerate}
\item We say that the $n$ species are \emph{stochastically persistent} if there exists a unique invariant probability measure $\Pi$ on $B_+ \times E$ such that $P_t((x,i),\cdot)$ converges in total variation to $\Pi$ for all initial conditions $(x,i) \in B_+\times E$.

\item Species $1,\ldots, k$ with $k<n$ are \emph{weakly persistent} if there is a unique invariant probability measure $\Pi^k$ on $B^k_+ \times E$ such that $P_t((x,i),\cdot)$ converges weakly to $\Pi^k$ for all initial conditions $(x,i) \in B_+\times E$.

\item Species $i$ \emph{goes extinct almost surely exponentially quickly with rate $\lambda<0$} if for all $(x,j) \in B_+ \times E$, 
$$P_{(x,j)} \left( \lim_{t\to \infty} \frac{\ln X_i^x(t)}{t} = \lambda\right)=1.$$
\end{enumerate}

We refer to \cite{Sch12} for an overview of the different forms of persistence and related notions.

Returning to our model, we note that intuitively due to the form of the $G^j$s, if species $i$ goes extinct, species $i+1$ must become extinct too since its only source of food is species $i$. Arguing in the same manner we conclude that species $i+2,\ldots,n$ must also become extinct. We shall see that this is indeed the case in the main results.

To do this, let us first recall that the average vector fields is defined by $G^\nu = \sum_{j=1}^N \nu_j G^j$ where $\nu = (\nu_1, \ldots, \nu_N)$ is the unique invariant probability measure of $J$ on $E$.
It is important to note that $G^\nu$ is also a Lotka-Volterra food chain defined by \eqref{functions_G^j}, that is $G_i^\nu(x) = x_i F_i^\nu(x)$ where $F^\nu$ is defined by \eqref{fucntions_F_i^j} with coefficients 
\begin{align}\label{a_ij_form_p_vector_field}
a_{ik}^\nu := \sum_{j=1}^N \nu_j \,a_{ik}^j.
\end{align} 
Remark also that $a_{ik}^\nu=a_{ik}$ for $k\neq 0$.

We say that a point $x$ is an {\it equilibrium} of a function $H$ provided $H(x) = 0$. Moreover we call it {\it positive} if all its coordinates are strictly positive (i.e. $x_i>0$ for all $i$).

\begin{lemma}
There exists a unique equilibrium for $F^\nu$.
\end{lemma}

This lemma has already been proved in \cite[Proposition 1]{LV_degenerate} where an explicit formula is given in the proof. We can also refer to \cite[Section 4]{HeningDang18b} where the authors also explicitly calculate the equilibrium in the case where $a_{ii}=0$ for $i=2,\ldots,n$.

For $k<n$, we let $F^\nu_{|k}$ denote the average vector field for the model \eqref{model} with $k$ species (i.e. $F^\nu$ defined by \eqref{G_nu} with $n$ replaced by $k$). This is important to remark that by the previous lemma, for each $1\leq k \leq n-1$, there exists a unique equilibrium for $F^\nu_{|k}$.

In fact, we will see that all the dynamics of the system will be determined by the equilibrium of  $F^\nu$ and the one of $F^\nu_{|k}$ for $k=1,\ldots,n-1$.
This is why for the remaining of the paper, we let $q^*$ denote the unique equilibrium of $F^\nu$ and $q^{*k}$ denote the one of $F^\nu_{|k}$ for $k=1,\ldots,n-1$. 
Note that these equilibria can have negative coordinates, i.e. $q^* \in \RR^n, q^{*k} \in \RR^k$.
 Next proposition investigates the relations between $q^{*k}$ and $q^{*k+1}$.

\begin{proposition}\label{properties_equil}
\begin{enumerate}
\item If $q^{*k+1}$ is positive then so too is $q^{*k}$.
\item The equilibrium $q^{*k}$ is positive if and only if $q^{*k}_k>0$.
\item If $q^{*k}_k \leq 0$ then $q^{*k+1}_{k+1} <0$. In particular if $q^{*k}$ is non-positive, then so too is $q^{*k+1}$.
\end{enumerate}
\end{proposition}

This proposition is a combination of Lemma 5.1 of \cite{LV_degenerate} and Lemma 4.1 of \cite{HeningDang18b}. 
Let us define the constant
$$\mathcal{I}^-_{k+1} : = -a_{k+1, 0}^\nu + a_{k+1,k}\, q^{*k}_k$$
which, as we shall see, is related to the survival of the species $k+1$.
It is important to remark that by Lemma 5.1 and the proof of Proposition 2.4 of \cite{LV_degenerate}, $\mathcal{I}^-_{k+1} = C \, q^{*k+1}_{k+1}$ for a constant $C>0$. Then by the second point of the previous proposition, $\mathcal{I}^-_{k+1}$ is positive if and only if $F^\nu_{|k+1}$ has a positive equilibrium.

Now, it is essential to notice that $F^\nu_{|1}$ always has a positive equilibrium since $a_{10}^j>0$ for all $j$. Then, we define $K$ as the first $k < n$ such that $F^\nu_{|K}$ has a positive equilibrium and $F^\nu_{|K+1}$ not or $K=n$ if there is no such $k$.

In the next theorems, we will see that there are three possibilities for the dynamic. The first one is when $K=n$ in which all species persist. The second one occurs under the conditions that $K=k<n$ and that $\mathcal{I}^-_{k+1} <0$ and we have the persistence of the $k$ first one and the extinction of the $n-k$ last species. The last possibility is the critical case when $K=k<n$ and $\mathcal{I}^-_{k+1} =0$ in which we have the persistence of the first $k$ species and the extinction of the last $n-k-1$ species while the $k+1$th dies in a more weaker sense.


\begin{theorem}[Persistence]\label{theorem_stoch_pers}
Suppose that Assumption \ref{Assumption_bio} holds and that $K=n$. Then,
\begin{itemize}
\item[(a)] There exists a unique invariant probability measure $\Pi$ on $B_+ \times E$.
\item[(b)] $\Pi$ is absolutely continuous with respect to the Lebesgue measure on $B \times E$ and for all $(x,j) \in B_+\times E$, $\Pi_t^{(x,j)}$ converges weakly to $\Pi$ a.s. as t goes to infinity.
\item[(c)] There exist a continuous function $W :B_+ \to \RR_+$ with $\lim\limits_{x \to \partial B} W(x) = \infty$ and constants $\rho,C>0$ such that for all $(x,j) \in B_+ \times E$
\begin{align*}
\|P_t((x,j),\cdot) -\Pi(\cdot) \|_{TV} \leq C (1 + W(x) )\, e^{-\rho t}.
\end{align*}

\end{itemize}
\end{theorem}

Further details on $W$ will be given in section \ref{Section_proof_thm_stoch_pers}. Note also that this theorem shows in particular the stochastic persistence of the $n$ species.

\begin{theorem}[Extinction]\label{theorem_extinction}
Suppose that Assumptions \ref{Assumption_bio} holds, that $K=k<n$ and that $\mathcal{I}^-_{k+1} <0$. Then, 
\begin{itemize}
\item[(i)] $\PP_{inv} (B \times E) = \PP_{inv} (B^k  \times E)$ and there exists a unique invariant probability measure $\Pi^k$ on $B_+^k \times E$. 
\item[(ii)] Species $k+1,k+2,\ldots,n$ go extinct almost surely exponentially quickly with rates $\mathcal{I}^-_{k+1}$, $-a^\nu_{k+2,0}, \ldots, -a^\nu_{n0}$ respectively.
\item[(iii)] For all $(x,j) \in B_+ \times E$, the transition kernel $P_t((x,j),\cdot)$ converges weakly to $\Pi^k(\cdot)$.
\end{itemize}

\end{theorem}

 In particular, species $1,\ldots, k$ are weakly persistent.

\begin{theorem}[Critical case]\label{theorem_critical_case}
Suppose that Assumptions \ref{Assumption_bio} holds, that $K=k<n$ and that $\mathcal{I}^-_{k+1} =0$. Then,
\begin{itemize}\item Species $1,\ldots,k$ are persistence in mean, i.e. $\lim_{t\to \infty} \frac{1}{t} \int_0^t X_i(s) ds = q_i^{*k}$ for $i=1,\ldots,k$ where $q^{*k}$ is the positive equilibrium of $F^\nu_{|k}$.
\item Species $k+2,k+3,\ldots,n$ go extinct almost surely exponentially quickly with rates $-a^\nu_{k+2,0},-a^\nu_{k+3,0}, \ldots, -a^\nu_{n0}$ respectively.
\item $\PP_{inv}(B_+^{k+1}\times E)$ is empty and the weak limit points of $\Pi_t$ lies almost surely in $\PP(B^k \times E)$.
\end{itemize}
\end{theorem}

It is important to notice that very few has been done for the critical case. We can cite \cite{NguyenStrickler2020} in which a criterion and numerous examples are
 given.

\section{Mathematical tools} \label{Sect_Math_tolls_stoch_pers}

In this section, we introducesome of the tools and properties needed for the proofs of the main theorems. We consider the general PDMP given by \eqref{model} assuming that the $F^j$s are $\mathcal{C}^\infty$ and sufficiently good for the process to be well defined at all $t\geq 0$, but are not necessarily of the form \eqref{fucntions_F_i^j}.


We define for $g:B\times E \to \RR$ smooth in the first variable the {\it infinitesimal generator} $\mathcal{L}$ by
%
$$\mathcal{L}g(x,i) = \langle G^i(x),\nabla g(x,i)\rangle + \sum_{j\in E} b_{ij}\left(g(x,j)-g(x,i)\right),$$
where $\langle\cdot,\cdot\rangle$ stands for the canonical scalar product of $\RR^n$ and $b_{ij}$s are the jump rates defined by \eqref{model} for $i\neq j$ and $b_{ii}=0$ otherwise.



\begin{definition}
A point $(y,i) \in B \times E$ is accessible from $(x,j) \in B \times E$ if for every neighbourhood $U$ of $(y,i)$, there exists $t \geq 0$ such that $P_t((x,j),U)>0$.

We denote by $\Gamma_{(x,j)}$ the set of points $(y,i)$ that are accessible from $(x,j)$ and for $D \subset B \times E$, we let $\Gamma_D = \underset{(x,j)\in D}{\bigcap} \Gamma_{(x,j)}$ be the set of accessible points from $D$.
\end{definition}

To characterize accessibility, we consider the deterministic control system associated to \eqref{model},
\begin{align} \label{equa_control}
\dot{y}(t) = \sum_{j=1}^N u^j(t) G^j(y(t)),
\end{align} 
where the control function $u = (u^1, \ldots, u^N) : \RR_+ \to \RR^N$ is at least piecewise continuous. We let $y(u,x,\cdot)$ denote the maximal solution to \eqref{equa_control} starting from $x$ with control function $u$.
The next characterization is Proposition 6.2 in \cite{B2018}.

\begin{proposition} \label{prop_control}
Suppose that $(x,j),(y,i) \in B\times E$. Then $(y,i) \in \Gamma_{(x,j)}$ if and only if for every neighbourhood $O \subset B$ of $y$, there exists a control function $u$ such that $y(u,x,t) \in O$ for some $t\geq 0$. 
\end{proposition}


We also need the so-called Hörmander conditions. Let 
\begin{align*}
\mathcal{F}_0 & := \{G^i-G^j \mid i,j = 1,\ldots,N\}, \\
\mathcal{F}_k & := \mathcal{F}_{k-1} \cup \left\{ [V,G^j] \mid V \in \mathcal{F}_{k-1}, \; j=1,\ldots,N \right\} \qquad k \ge 1,
\end{align*}
where $[\cdot,\cdot]$ denotes the Lie bracket operator. 
Recall that this operator is defined for every smooth vector fields $V,W : \RR^n \to \RR^n$ and every $x \in \RR^n$ by
$$[V,W](x) := DW(x) V(x)-DV(x) W(x),$$
where $DV(x) = \left(\frac{\partial V_i}{\partial x_j}(x)\right)_{i,j}$ stands for the Jacobian matrix of $V$ at $x$.

For $x \in \RR^n$ and $k \geq 1$, we write $\mathcal{F}_k(x) := \left\{ V(x) \mid V \in \mathcal{F}_k \right\}$.

\begin{definition}\label{def_hormander}
We say that $\eqref{model}$ satisfies the strong Hörmander condition at $(x,j) \in B \times E$ if there exists $k \in \NN$ such that 
$$span(\mathcal{F}_k (x)) = \RR^n.$$

It satisfies the weak Hörmander condition if $\mathcal{F}_0$ is weakened to $\mathcal{F}'_0 = \{G^1,\ldots, G^N\}$ and there exist $k$ and $(x,j)$ like in the strong Hörmander condition.
\end{definition}


By the form of the $G^j$s, spaces of the form 
$$B_+^I := \left\{ x \in B \mid x_i >0 \text{ if } i \in I \text{ and } x_i=0 \text{ otherwise} \right\}$$
for some $I \subset \left\{1,\ldots,n\right\}$ are invariant for the process given by \eqref{model}. An ergodic probability measure on the boundary must then be supported by such a space. More precisely for each $\mu \in \PP_{erg}(\partial B \times E)$, there exists $I_\mu \subset \left\{1,\ldots,n\right\}$ such that $\mu(B_+^{I_\mu} \times E) =1$.

For each $\mu \in \PP_{inv}(B \times E)$ and each $j \in E$, we define the measure $\mu^j$ by $\mu^j(A) := \mu(A\times \{j\})$ for all Borel set $A \subset B$. 
As an example consider the invariant measure $\mu$ defined by  $\mu^j= \nu_j \delta_0 \Ind_{\{j\}}$ for all $j$ where $\delta_0$ is the Dirac measure at $0 \in B$ and $\nu$ is the unique invariant probability measure on $E$. Then $\mu = \delta_0 \otimes \nu$ and $I_\mu = \emptyset$.

\begin{definition}\label{def_invasion_rate}
The invasion rate of species $i$ with respect to $\mu \in \mathcal{P}_{inv}(\partial B \times E)$ is defined by
\begin{align} \label{Invasion_rate_equ}
\lambda_i(\mu)  := \sum_{j=1}^N \int_B F_i^j(x) d\mu^j(x).
\end{align}
%
\end{definition}

The next proposition gives a very useful property of the invasion rates with respect to an ergodic probability measure on the boundary.


\begin{proposition}\label{Invasion_rate_0}
Let $\mu \in \PP_{erg}(\partial B \times E)$. Then for $i\in I_\mu$
$$\lambda_i(\mu) =0.$$
\end{proposition}

\prf
Let the process $Z^\mu(t) = (X^\mu(t), J^\mu(t))$ be the solution of \eqref{model} with initial distribution $\mu \in \PP_{erg}(\partial B \times E)$. Since $F_i^j$s are $\mathcal{C}^\infty$ and $B \times E$ is compact, then $F_i^j$ are $\mu$-integrable. For $i\in \left\{ 1,\ldots,n\right\}$ by Birkhoff's ergodic theorem we get that
%
%
\begin{align*}
\lim_{t \to \infty} \frac{\ln \left( X_i^\mu(t)  \right)}{t} & =  \lim_{t \to \infty} \frac{1}{t} \int_0^t  \frac{\dot X_i^\mu (s)}{ X_i^\mu (s)} \, ds = \lim_{t \to \infty} \frac{1}{t} \int_0^t F_i^{J^\mu(s)}(X^\mu(s)) \, ds\\
&  = \lim_{t \to \infty} \frac{1}{t} \int_0^t \sum_{j=1}^N \Ind_{\{J^\mu(s) = j\}} F_i^j (X^\mu(s)) \, ds \\
& =  \sum_{j=1}^N \int_{B } F_i^j (x) \, \mu^j(dx) = \lambda_i(\mu).
\end{align*}
For $i \in I_\mu$, if $\lambda_i(\mu) \neq 0$ the previous equation shows that 
$$\lim_{t\to\infty} X_i^\mu(t) = \left\{ \begin{array}{ll}
	0		& \text{ if } \lambda_i(\mu) <0,\\
	\infty	& \text{ if } \lambda_i(\mu) >0,
    \end{array}
	\right.
$$
which is in contradiction with $\mu(B_+^\mu \times E) =1$.

{\hfill$\square$\bigbreak}

\begin{remark}
This proof doesn't require that $F^j$ is $\mathcal{C}^\infty$ but only $\mathcal{C}^0$, thus the previous result is much more general.

It is also interesting to note that this result can be deduced from the equality $\mu \mathcal{L} g = 0$ with $g(x,i) = \ln(x_i)$.

Note again that this has already been proved in the setting of discrete Markov chains in \cite{SBA11} and in the context of some general stochastic differential equations in \cite{B2018,hening2020general}, but not in the context of switching vector fields.

\end{remark}

The next proposition shows a link between this invasion rate and the empirical occupation measure.

\begin{proposition}\label{ln_sur_t_cv_lambda_pi_t}
For all $(x^*,j^*) \in B_+ \times E$, 
\begin{align}
\P_{(x^*,\,j^*)}\left( \lim\limits_{t \to \infty} \left( \frac{\ln(X_i(t))}{t} -\lambda_i(\Pi_t) \right) =0 \right) =1.
\end{align}

\end{proposition}

\prf
Like in Proposition \ref{Invasion_rate_0}, we have that
\begin{align*}
\frac{\ln \left( X_i^{x^*}(t)  \right) - \ln \left(x^*  \right)}{t} &  =  \sum_{j=1}^N \frac{1}{t} \int_0^t \Ind_{\{J^{j^*}(s) = j\}} \, F_i^j \left(X^{x^*}(s)\right) \, ds \\
& = \sum_{j=1}^N \int_B F_i^j(x)\, \frac{1}{t} \int_0^t \Ind_{\{Z^{(x^*,\,j^*)}(s) \in dx\}}  \, ds \\
& =  \sum_{j=1}^N \int_B  F_i^j (x) \, \Pi^{(x^*,j^*)}_t(dx) = \lambda_i\left(\Pi^{(x^* ,\,j^*)}_t\right).
\end{align*}
We conclude the proof by letting $t$ go to infinity.

{\hfill$\square$\bigbreak}


When the $F^j$s are of the form \eqref{fucntions_F_i^j}, we can determine exactly the form of the set $I_\mu$ for $\mu \in \PP_{erg}(\partial B \times E)$. The next proposition has already been proved in \cite[Lemma A.1]{HeningDang18c} in the stochastic differential equations setting. We extend here the result to the model with random switching.

\begin{proposition}\label{support_ergodic}
Suppose that the $F^j$s are of the form \eqref{fucntions_F_i^j}, then for each $\mu \in \PP_{erg}(\partial B \times E)$ there exists $1\leq k < N$ such that $I_\mu = \left\{1,2,\ldots,k\right\}$ or $I_\mu = \emptyset$.
\end{proposition}

\prf 
We start by noting that $I_\mu  =\emptyset$ if and only if 
$\mu = \delta_0 \otimes \nu$.
We then suppose that $\mu \neq \delta_0 \otimes \nu$ and we write $I_\mu = \left\{i_1,\ldots,i_k\right\}$ with $i_1 <i_2 < \ldots < i_k$. We first assume that $i_1 >1$, then by definition $\mu \left(\{(x,j) \in B\times E \mid x_{i_1-1} =0\} \right) = 1$. By Proposition \ref{Invasion_rate_0}, 
\begin{align*}
0 = \lambda_{i_1} (\mu) & =   - a_{i_1,0}^\nu - a_{i_1,i_1} \sum_{j=1}^N  \int x_{i_1} d\mu^j(x) - a_{i_1,i_1 +1 } \sum_{j=1}^N  \int x_{i_1+1} d\mu^j(x) 
\end{align*}
where $a_{i_1,0}^\nu$ is defined by \eqref{a_ij_form_p_vector_field}. Since the right-hand side is strictly negative, this is a contradiction. We now assume that there exists $m$ such that $i_m -1> i_{m-1}$. Then $\mu \left(\{(x,l) \in B\times E \mid x_{i_m-1} =0\} \right) = 1$ and once again by Proposition \ref{Invasion_rate_0},
\begin{align*}
0 = \lambda_{i_m} (\mu) & = -a_{i_m,0}^\nu - a_{i_m,i_m} \sum_{j=1}^N \int x_{i_m} d\mu^j(x) - a_{i_m,i_m +1} \sum_{j=1}^N \int x_{i_m+1} d\mu^j(x) 
\end{align*}
which also leads to a contradiction since the right-hand side is always strictly negative.

{\hfill$\square$\bigbreak}

\section{Proofs}\label{Section proof}
			
\subsection{Proof of Theorem \ref{theorem_stoch_pers}} \label{Section_proof_thm_stoch_pers}

Theorem \ref{theorem_stoch_pers} will be deduced from the next one which is Corollary 6.3 of \cite{B2018}.

\begin{theorem} \label{Criter_thm_stoc_persi}
Suppose that there exist $\alpha_1,\ldots,\alpha_n >0$ such that for all $\mu \in \PP_{erg}(\partial B \times E)$
\begin{align}\label{Criter_alpha_i_lambda_i}
\sum_{i=1}^n \alpha_i\lambda_i(\mu) >0.
\end{align}
Assume moreover that there exists $(p^*,j) \in  \Gamma_{(B_+ \times E)} \cap (B_+ \times E)$ which satisfies the strong Hörmander condition. Then 
\begin{itemize}
\ita There exists a unique invariant probability measure $\Pi$ on $B_+ \times E$.

\itb $\Pi$ is absolutely continuous with respect to the Lebesgue measure of $B \times E$ and for all $(x,j) \in B_+ \times E$, $\Pi_t^{(x,j)}$ converges weakly to $\Pi$ a.s. as t goes to infinity.

\itc For all $(x,j) \in B_+ \times E$,
$$ \| P_t((x,j), \cdot) -\Pi(\cdot)\|_{TV} \leq C (1+W(x)) e^{-\rho t} $$
for some $\rho>0$ and 
$$W(x) :=  e^{( \theta \max \{-\sum_{i=1}^n \alpha_i \ln(x_i),1\})}$$
for some $\theta>0$.

\end{itemize}
\end{theorem}

%



As said in the introduction, the positive equilibrium of $F^\nu$, $q^*$, will play a central role and more precisely its variations $(q^*,j)$, $j \in E$ are points that will satisfy hypotheses of the previous theorem. We start by proving the accessibility of these points.


\begin{proposition}\label{Accessibiliy_q*}
Assume that $F^\nu$ has a unique positive equilibrium $q^*$. Then for any $j\in E$, $(q^*,j)$ is accessible for all $(x,i) \in B_+ \times E$.
\end{proposition}

\prf Since $G^\nu$ is still a Lotka-Volterra food chain, by Theorem 5.3.1 of \cite{HS98}, $q^*$ is globally asymptotically stable for the dynamic induced by the system 
$$\dot Y_i(t)  = G_i^\nu(Y(t)) \qquad i=1,\ldots,n.$$ 
This means that for every initial condition $y_0 \in B_+$, the solution $Y^{y_0}(t)$ has the property that $Y^{y_0}(t) \xrightarrow[t \to \infty]{} q^*$.

Then by Proposition \ref{prop_control} it is enough to choose $u(t) = \nu $ which is the invariant probability measure on $E$ 
for the control function in \eqref{equa_control} to have accessibility of $(q^*,j)$ for all $j \in E$.

{\hfill$\square$\bigbreak} 

\begin{remark}
This proof also shows that if a point $p^*$ is globally asymptotically stable for the dynamic induced by a linear combination of the vector fields $G^j$, then for all $k \in E$, $(p^*,k)$ is accessible for all $(x,i) \in B_+\times E$.
\end{remark}


We go on with the strong Hörmander condition.

\begin{proposition} \label{Hormander_at_q*}
Under assumptions of Theorem \ref{theorem_stoch_pers}, the strong Hörmander condition  
holds true at $(q^*,j)$ for all $j\in E$.
\end{proposition}

\prf 
By Assumption \ref{Assumption_bio}, there exist $\beta_1,\beta_2$ such that $G_i^{\beta_1} = G_i^{\beta_2}$ for $i=2,\ldots,n$ and $G_1^{\beta_1}(q^*)-G_1^{\beta_2}(q^*) \neq 0$. Let us consider the sequence $(b^k)$ defined by
$$\left\{\begin{array}{lcl}
b^1 = G^{\beta_1}-G^{\beta_2} = (G_1^{\beta_1}-G_1^{\beta_2}) e_1, \\
b^{k+1} = [b^k,G^{\beta_1}]
\end{array} \right.$$
where $e_1$ is the first vector of the canonical basis of $\RR^n$. We want to prove that $\{b^1(q^*), \ldots, b^n(q^*)\}$ is a basis of $\RR^n$.

To do this, we first remark that $DG^{\beta_1}$ is a tridiagonal matrix
$$DG^{\beta_1} = \begin{pmatrix}
g_{11} & g_{12}  & 0 & \cdots & 0 \\
g_{21}  & g_{22} & g_{23} & \ddots & \vdots \\
0 & g_{32} & \ddots & \ddots & 0 \\
\vdots & \ddots & \ddots & \ddots & g_{n-1,n}\\
0 & \cdots & 0 & g_{n,n-1} & g_{nn}
\end{pmatrix}$$
with coefficients
$$\left\{\begin{array}{ll}
g_{ii}(x) = F^{\beta_1}_i (x) - a_{ii}\, x_i,\\
g_{i,i+1}(x) = - a_{i,i+1} x_i,\\
g_{i,i-1}(x) = a_{i,i-1} x_i.
\end{array}\right.
$$
It is then not hard to see that for each $x \in B_+$,
$$ b^k_k(x) = \left( G^{\beta_1}_1 (x)-G^{\beta_2}_1(x) \right) \, \prod_{i=2}^k a_{i,i-1} x_i \quad \text{ and } \quad b^k_i (x) = 0 \text{ for } i >k.
$$
By Assumption \ref{Assumption_bio}, $\{b^1(q^*), \ldots, b^n(q^*)\}$ is a basis of $\RR^n$ and therefore the strong Hörmander condition holds true at $(q^*,j)$ for all $j \in E$.
%

{\hfill$\square$\bigbreak}

\begin{remark}
If Assumption \ref{Assumption_bio} is replaced by
$$\left\{ \begin{array}{ll}
a_{n0}^{\beta_1} \neq a_{n0}^{\beta_2}, \\
a_{i0}^{\beta_1} = a_{i0}^{\beta_2} & \text{ for } i=1,\ldots,n-1.
\end{array} \right.$$
Then a very similar proof also shows the strong Hörmander condition at $(q^*,j)$ for all $j\in E$.
\end{remark}

\begin{remark}\label{rem_H_faible}
This proof works also in a more general context but gives the weak Hörmander condition. Indeed, suppose that $G^1$ is a food chain living on $\RR_+^n$ having the properties that species $i$ eats species $i-1$, is hunted by species $i+1$, can have intra-specific interactions and doesn't interact with the other species. In a mathematical point of view, this means that for all $x$ in  the interior of $\RR_+^n$,
\begin{align*}
\left\{  
\begin{array}{ll}
\frac{\partial G^1_i}{\partial x_{i+1}}(x) <0 & \text{ for } i= 1,\ldots,n-1, \\
\frac{\partial G^1_i}{\partial x_{i-1}}(x) >0 & \text{ for } i= 2,\ldots,n, \\
\frac{\partial G^1_i}{\partial x_j}(x) =0 & \text{ for } j \neq i-1,i,i+1. 
\end{array}
\right.
\end{align*}
We consider a second vector fields 
$$G^2(x) := h(x) e_1$$ 
where $h : \RR^n \to \RR$ is a smooth function whose induced dynamics leaves $\RR_+$ invariant (i.e. $h(0)=0$) and $e_1$ stands for the first vector of the canonical basis of $\RR^n$. We consider the randomly switched model given by \eqref{model} for the vector fields $G^1$ and $G^2$. 
Then the weak Hörmander condition holds true at $(x,j)\in \RR_+^n \times \{1,2\}$ provided $h(x) \neq 0$.

Indeed, let's take the same family of vector fields $(b^k)$ as in the previous proposition, i.e.
$$\left\{\begin{array}{ll}
b^1  = G^2, \\
b^{k+1} = [b^k,G^1],
\end{array} \right.$$
which has also the property that
$$
b^k_k(x) = h(x) \prod_{i=2}^k \frac{\partial G^1_i}{\partial x_{i-1}}(x) \quad \text{ and } \quad b^k_i (x) = 0 \text{ for } i >k.
$$
Once again, this makes the set $\{b^1(x), \ldots, b^n(x)\}$ a basis of $\RR^n$ provided $h(x) \neq 0$ and proves the weak Hörmander condition at $(x,j)$.
\end{remark}

\begin{remark}
Under the assumptions of Theorem \ref{theorem_stoch_pers} and with Remark \ref{rem_H_faible}, by using Proposition  $2.10$ of \cite{BHS18} we also find the strong Hörmander condition at $(q^*,j)$.


\end{remark}


The next result is a key one to characterize the persistent condition \eqref{Criter_alpha_i_lambda_i}.

\begin{proposition} \label{delta_pos_iif_existe_p_i_st_sum_pi_mu_lmada_i_pos}
Condition \eqref{Criter_alpha_i_lambda_i} of Theorem \ref{Criter_thm_stoc_persi} is equivalent to the existence of a unique positive equilibrium for $F^\nu$.
\end{proposition}

\prf
We first suppose the existence of a positive equilibrium for $F^\nu$.
By Lemma 4 in \cite{SBA11}, condition \eqref{Criter_alpha_i_lambda_i} is equivalent to 
\begin{align}\label{lambda_mu}
\max_{i=1, \ldots, n} \lambda_i(\mu) >0, 
\end{align}
for every $\mu \in \PP_{inv}(\partial B \times E)$. We first prove \eqref{lambda_mu} for
$\mu \in \PP_{erg}(\partial B \times E)$.

If $\mu = \delta_0 \otimes \nu$, then $\lambda_1(\mu) = a_{10}^\nu >0$. If $\mu$ is any other ergodic probability measure, let $\overline{x}_i =\sum_{j=1}^N \int x_i \mu^j(dx)$.
By Proposition \ref{support_ergodic}, there exists $1\leq k<n$ such that $\mu(B_+^k \times E) =1$.  Then,
\begin{align*}
\lambda_{k+1}(\mu) & = -a_{k+1,0}^\nu + a_{k+1,k}\,\overline{x}_k\\
& = -a_{k+1,0}^\nu + a_{k+1,k} \, q^{*k}_k
\end{align*}
where we recall that $q^{*k}$ is the equilibrium of $F^\nu_{|k}$. With Proposition \ref{Invasion_rate_0}, we obtain the system of equation $\lambda_i(\mu) = 0$ for $i=1, \ldots,k$. The solution of this system is exactly the equilibrium of $F^\nu_{|k}$, then $\overline{x}_k =q^{*k}_k$ and the second equality is therefore proved.

By Lemma 5.1 and the proof of Proposition 2.4 in \cite{LV_degenerate}, $\lambda_{k+1}(\mu) = C \, q^{*k+1}_{k+1}$ for some constant $C>0$. Then the existence of a positive equilibrium for $F^\nu$ and Proposition \ref{properties_equil} prove the positivity of $\lambda_{k+1}(\mu)$.

%

Suppose now that $\mu \in \PP_{inv}(\partial B \times E)$. Then by the ergodic decomposition theorem there exist $\mu_{i_1},\ldots,\mu_{i_k} \in \PP_{erg}(\partial B \times E)$ and $p_j >0$ such that
$$\mu = p_1 \mu_{i_1} + \ldots + p_k \mu_{i_k}$$
with $i_1<\ldots<i_k$, $\mu_{i_j} (B_+^{i_j} \times E) =1$ and $\sum_{j=1}^k p_j =1$. Thus by Proposition \ref{Invasion_rate_0}, $\lambda_{i_1+1}(\mu_{i_j}) = 0$ for $j=2,\ldots,k$. By the ergodic case, this implies that $\lambda_{i_1+1}(\mu) = p_1 \lambda_{i_1+1}(\mu_{i_1}) >0$.

Conversely, suppose that condition \eqref{Criter_alpha_i_lambda_i} holds. Then by Theorem 6.1 and Theorem 4.4 of \cite{B2018}, there exists at least one ergodic probability measure $\mu$ supported by $B_+\times E$. By mimicking the proof of Proposition \ref{Invasion_rate_0}, we get that 
$$\lambda_i(\mu) =0  \qquad i=1,\ldots,n.$$
Let $\overline{x}_i =\sum_{j=1}^N \int x_i \mu^j(dx)$, then the previous system can be rewritten as follows
$$\left\{
     \begin{array}{ll}
     0 = a_{10}^\nu -a_{11} \overline{x}_1 - a_{12} \,\overline{x}_2 								& i=1,\\
     0 = -a_{i0}^\nu + a_{i,i-1} \overline{x}_{i-1} - a_{ii} \, \overline{x}_i - a_{i,i+1} \overline{x}_{i+1} 	& i= 2, \ldots,n-1, \\
     0 = -a_{n0}^\nu + a_{n,n-1} \overline{x}_{n-1} -a_{nn} \,\overline{x}_n 					& i=n.
     \end{array}
\right. $$
This system of equations has exactly one solution $(\overline{x}_1,\ldots,\overline{x}_n)$ which is by linearity the equilibrium of $F^\nu$. 
Moreover since the support of $\mu$ is $B_+ \times E$, $\overline{x}_i>0$ for all $i$ and this concludes the proof.

{\hfill$\square$\bigbreak}

The existence of a unique positive equilibrium for $F^\nu$ implies condition \eqref{Criter_alpha_i_lambda_i} by Proposition \ref{delta_pos_iif_existe_p_i_st_sum_pi_mu_lmada_i_pos}. It also implies accessibility of $(q^*,j)$ for all $j\in E$ by Proposition \ref{Accessibiliy_q*}. Moreover by Assumption \ref{Assumption_bio} and Proposition \ref{Hormander_at_q*}, the strong Hörmander condition holds true at these points. Then Theorem \ref{theorem_stoch_pers} follows from Theorem \ref{Criter_thm_stoc_persi}.


\subsection{Proof of Theorem \ref{theorem_extinction}}\label{sec_extinction}

\paragraph{Part i)}


Under the hypotheses of Theorem \ref{theorem_extinction}, there is no invariant probability measure on $B^i_+ \times E$ for $i=k+1, \ldots,n$.
Indeed, let us assume that there is one supported by $B^i_+ \times E$ for one $i=k+1,\ldots,n$. By mimicking the proof of Proposition \ref{delta_pos_iif_existe_p_i_st_sum_pi_mu_lmada_i_pos}, we can show that $F^\nu_{|i}$ has a unique positive equilibrium so that Proposition \ref{properties_equil} gives a contradiction.

We obtain in particular that
$$\PP_{inv}(B \times E) = \PP_{inv}(B^k \times E)$$
which concludes the proof of the first part.

\paragraph{Part ii)} 
%
%
%


We start by claiming that for all $\mu \in \PP_{erg}(B^k \times E)$ 
\begin{align}\label{calcul_lmaba_i_extinction}
\lambda_i(\mu) = \left\{\begin{array}{ll}
-a^\nu_{k+1,0} & \text{ if } I_\mu=\{1,\ldots,j\}, \; j<k \text{ and } i=k+1, \\
 \mathcal{I}^-_{k+1} &  \text{ if }  I_\mu=\{1,\ldots,k\} \text{ and } i=k+1, \\
-a^\nu_{i0} & \text{ if } i>k+1.
\end{array} \right.
\end{align}
By definition of $\lambda_i$, bottom and upper cases are immediate. Let $\mu$ be such that $I_\mu=\{1,\ldots,k\}$, then by Proposition \ref{Invasion_rate_0}, $\lambda_i(\mu) =0$ for $i=1,\ldots,k$. The solution of this system of equations is $q^{*k}$ and then $\lambda_{k+1}(\mu)  = -a_{k+1, 0}^\nu + a_{k+1,k}\, q^{*k}_k = \mathcal{I}^-_{k+1}.$ The claim is thus proved. We recall that by the assumptions, $\mathcal{I}^-_{k+1}<0$.

It follows from \eqref{calcul_lmaba_i_extinction} and Proposition \ref{ln_sur_t_cv_lambda_pi_t} that for all $(x,j) \in B_+ \times E$,
$$\lim\limits_{t\to \infty} X^x_i(t) = 0 \qquad i=k+1,\ldots,n.$$

Now we claim that with probability $1$ the weak limit points of $(\Pi_t^{(x,j)})$ when $t\to \infty$ are in $\PP_{inv}(B_+^k \times E)$ for all $(x,j) \in B_+ \times E$.

Indeed by Theorem 2.2 of \cite{B2018}, for every $(x,j) \in B_+ \times E$, $(\Pi_t^{(x,j)})$ is tight and the weak limit points are invariant probability measure almost surely. Suppose by contradiction that there is a point $(x,j) \in B_+ \times E$ such that with positive probability, there exists a sequence $(t_m)_m$ such that $\Pi^{(x,j)}_{t_m} \to \mu$ weakly, $\mu = \rho_1 \mu_1 + \rho_2 \mu_2$ with $\mu_1 \in \PP_{inv}(B^{k-1} \times E)$, $\mu_2 \in \PP_{inv} (B_+^k \times E)$, and $\rho_1>0$.

In the same manner as in the proof of Proposition \ref{delta_pos_iif_existe_p_i_st_sum_pi_mu_lmada_i_pos}, $\lambda_k(\mu_1) >0$.  Moreover Proposition \ref{Invasion_rate_0} implies that $\lambda_k(\mu_2) = 0$. Then with Proposition \ref{ln_sur_t_cv_lambda_pi_t} we get that
$$\lim_{m\to\infty} \frac{\ln X^x_k(t_m)}{t_m} = \lim_{m\to \infty} \lambda_k(\Pi^x_{t_m}) = \lambda_k(\mu) >0.$$
This is in contradiction with the compactness of $B$ and the claim is therefore proved.

By Proposition \ref{ln_sur_t_cv_lambda_pi_t}, equation \eqref{calcul_lmaba_i_extinction} and the property of $(\Pi_t^{(x,j)})$ above we finally get almost surely that
$$\lim_{t\to \infty} \frac{\ln X_i^x(t)}{t} = \left\{ \begin{array}{ll}
-a^\nu_{i0} & \text{ if } i>k+1, \\
\mathcal{I}^-_{k+1} & \text{ if } i=k+1.
\end{array}\right.$$
Then species $k+1,k+2,\ldots,n$ go extinct almost surely exponentially quickly with rate $\mathcal{I}^-_{k+1}$, $-a^\nu_{k+2,0}, \ldots, -a^\nu_{n0}$ respectively.

\paragraph{Part iii)} We start this part by the next proposition which gives a characterisation of the remaining species $x_1,\ldots, x_k$.

\begin{proposition} \label{persitence_en_proba}

For every $\varepsilon >0$, there exists a compact set $K \subset B_+^k \times E$ such that for all $(x,j) \in B_+ \times E$,
$$\liminf_{t\to \infty} \P_{(x,j)}\left( (X_1(t),\ldots,X_k(t),J(t)) \in K \right) \geq 1- \varepsilon.$$
\end{proposition}

\begin{remark}
This property is called persistence in probability of species $1,\ldots,k$ in \cite{HeningDang18c} and we still refer to \cite{Sch12} for an overview of the different forms of persistence.

\end{remark}


\prf
Since $F^\nu_{|k}$ has a unique positive equilibrium, by Proposition \ref{delta_pos_iif_existe_p_i_st_sum_pi_mu_lmada_i_pos} there exist $\alpha_1,\ldots,\alpha_k >0$ such that
$$\sum_{i=1}^k \alpha_i\lambda_i(\mu) >0$$ 
for every $\mu\in \PP_{inv}(\partial B^k\times E)$.
Let $\xi : \RR \to \RR_+$ be a smooth function with bounded first and second derivatives such that $\xi(t) = t$  for $t \geq 1$ and let
$$V(x,j) := \xi\left(-\sum_{i=1}^k \alpha_i \ln(x_i)\right)\!.$$
The latter with the help of Lemma $8.2$ and Theorem $6.1$ of \cite{B2018} (see also the proof of Theorem 1.14 of \cite{StricklerPhD} or Lemma 3.15 of \cite{BL16}) implies that there exist constants $T,\kappa,\theta>0$ and $0<\rho<1$ such that the function 
$$V_\theta := e^{\theta V}$$
satisfies on $B_+ \times E$
\begin{align*}
P_T V_\theta \leq \rho V_\theta + \kappa.
\end{align*}
Let $(x,j) \in B_+ \times E$, then by the strong Markov property
\begin{align*}
\E_{(x,j)}\left[V_\theta( Z(mT))\right] \leq \rho^m V_\theta(x,j) +\kappa \sum_{i=0}^{m-1} \rho^i,
\end{align*}
and therefore
\begin{align}\label{equa_pers_in_proba_1}
\limsup_{m\to \infty} \E_{(x,j)}\left[V_\theta( Z(mT))\right] \leq \frac{\kappa}{1-\rho}.
\end{align}
It is straightforward that
\begin{align*}
\mathcal{L} V_\theta(x,j) = -\theta V_\theta(x,j) \sum_{i=1}^k \alpha_i F_i^j(x).
\end{align*}
By compactness of $M$ and smoothness of the $F^j$s, there exists $c>0$ such that
\begin{align*}
\mathcal{L} V_\theta \leq c V_\theta.
\end{align*}
For $(x,j)$ fixed, let the stopping time $\tau_m := \inf \left\{ t \geq 0 \mid V_\theta(Z^{(x,j)}(t)) \geq m \right\}$. For $a,b$, we write $a \wedge b$ for the minimum between $a$ and $b$. Then by Dynkin's formula
\begin{align*}
\E_{(x,j)} \left[V_\theta (Z(t \wedge \tau_m)) \right] & =V_\theta(x,j) + \E_{(x,j)} \left[ \int_0^{t\wedge \tau_m} \mathcal{L}V_\theta (Z(s)) ds\right] \\
& \leq V_\theta(x,j) + c \,\E_{(x,j)} \left[ \int_0^{t\wedge \tau_m} V_\theta (Z(s)) ds\right].
\end{align*}
Since $V_\theta(x,j) \to \infty$ when $x \to \partial B^k$, by letting $m \to \infty$, we obtain that
\begin{align*}
\E_{(x,j)} \left[V_\theta (Z(t)) \right] & \leq V_\theta(x,j) + c \,\E_{(x,j)} \left[ \int_0^{t} V_\theta (Z(s)) ds\right],
\end{align*}
and by using Gronwall's lemma, we get that
\begin{align*}
\E_{(x,j)} \left[V_\theta (Z(t)) \right] \leq e^{ct}\, V_\theta(x,j).
\end{align*}
For $t \in [mT,(m+1) T]$, the strong Markov property implies that
\begin{align}\label{equa_pers_in_proba_2}
\E_{(x,j)} \left[V_\theta (Z(t)) \right] \leq e^{cT}\, \E_{(x,j)} \left[V_\theta (Z(mT)) \right].
\end{align}
Finally, we combine \eqref{equa_pers_in_proba_1} and \eqref{equa_pers_in_proba_2} to obtain
\begin{align}\label{equa_pers_in_proba_3}
\limsup_{t\to\infty} \E_{(x,j)} \left[V_\theta (Z(t)) \right] \leq  e^{cT}\frac{\kappa}{1-\rho}.
\end{align}
Now, let us fix $\varepsilon>0$ and define the compact set 
$$K := \left\{  (x,j) \in B^k_+ \times E \mid \varepsilon \, V_\theta(x,j) \leq  e^{cT}\frac{\kappa}{1-\rho} \right\}.$$
Then with \eqref{equa_pers_in_proba_3}, we get that
\begin{align*}
\limsup_{t\to \infty} \P_{(x,j)}\left(Z(t) \not \in K\right) \leq \limsup_{t\to \infty} \E_{(x,j)}\left[ \varepsilon\, e^{-cT} \frac{1- \rho}{\kappa} \, V_\theta(Z(t))\right] \leq \varepsilon,
\end{align*}
which concludes the proof.

{\hfill$\square$\bigbreak}

By hypotheses of Theorem \ref{theorem_extinction}, $F^\nu_{|k}$ has a unique positive equilibrium which implies by Theorem \ref{theorem_stoch_pers} that there exists a unique invariant probability measure $\Pi^k$ on $B_+^k \times E$ such that the semi-group converges exponentially quickly to it. Moreover, with the form of the function $W$ the convergence is uniform on each compact set of $B_+^k \times E$ meaning that for each continuous bounded function $f : B_+^k \times E \to \RR_+$ and each compact set $K\subset B_+^k \times E$,
\begin{align}\label{cv_unif_compact}
\lim_{t\to \infty}\left(\sup_{(x,j)\in K} \left|P_tf(x,j)-\Pi^k f(x,j)\right|\right) = 0.
\end{align}
Moreover, by \cite[Proposition 2.1]{BMZIHP} the process is Feller, meaning that for each $f : B \times E \to \RR$ continuous bounded, the map $(t,(x,j)) \mapsto P_t f(x,j)$ is continuous.

Then, with the extinction exponentially quickly of species $k+1,\ldots,n$, Proposition \ref{persitence_en_proba}, \eqref{cv_unif_compact} and the Feller property, we can imitate the proof of Theorem $(iii)$ of \cite{HeningDang18c} to get the weak convergence of the semi-group $P_t((x,j), \cdot)$ to the invariant probability measure $\Pi^k(\cdot)$. This concludes the proof of the third part.


\subsection{Proof of Theorem \ref{theorem_critical_case}} 

The two first parts of Theorem \ref{theorem_critical_case} work in the same way as in the proof of Theorem \ref{theorem_extinction} combined with Proposition \ref{ln_sur_t_cv_lambda_pi_t}. We also refer to the proof of Theorem 1.1.(ii) of \cite{HeningDang18c}.

For the last part, we want to use Corollary 2.8 of \cite{NguyenStrickler2020}. To do so, let us define the extinction set 
$$M_0 := \left\{ (x,i) \in B\times E \mid x_{k+1}=0\right\},$$
the living set $M_+ := (B\times E) \backslash M_0$ and the map 
$$\begin{array}{ccccl}
 H & : & B \times E &\longrightarrow & \RR\\
 &&(x,i)& \longmapsto &H(x,i):= - F_{k+1}^i(x).
\end{array}$$
Then by Theorem 6.1 of \cite{B2018}, $H$, $M_0$ and $M_+$ satisfy the assumptions of \cite{NguyenStrickler2020} for some good function $V$. It thus remains to prove the following: \emph{if $\PP_{inv}(M_+)$ is not empty, then there is $\mu \in \PP_{inv}(M_+)$ and $\tilde\mu \in \PP_{inv}(M_0)$ such that} 
$$\mu H > \tilde\mu H.$$
First remark that with Proposition \ref{support_ergodic}, $\PP_{inv}(M_0) = \PP_{inv}(B^k \times E)$. Moreover, since $F^\nu_{|k}$ has a positive equilibrium and with the help of Theorem \ref{theorem_stoch_pers}, there is an ergodic probability measure $\tilde\mu$ such that $\tilde\mu(B_+^k \times E) =1$. Then by \eqref{calcul_lmaba_i_extinction}, $\lambda_{k+1}(\tilde\mu) = \mathcal{I}^-_{k+1}$.

On the other hand, if there is $\mu \in \PP_{inv}(M_+)$, then we can always choose it ergodic and with Proposition \ref{Invasion_rate_0} and Proposition \ref{support_ergodic}, $\mu(B_+^{k+1} \times E) >0$. Then, the same argument as in Proposition \ref{delta_pos_iif_existe_p_i_st_sum_pi_mu_lmada_i_pos} implies that $\mathcal{I}^-_{k+1} >0$ and thus $\tilde\mu H <0$. 

We conclude the proof by noting that with the help of Lemma 7.5 of \cite{B2018} or of the proof of Proposition 2.7 of \cite{NguyenStrickler2020}, $\mu H =0$ for all $\mu \in \PP_{inv}(M_+)$.


\subsection{Proof of Proposition \ref{compact_B}} \label{proof_compact_B}

%
%
%

Let $$S(t) := X_1(t) + \sum_{i=2}^n \varepsilon_i X_i(t)$$
be the total weighted sum of all individuals of the population with weights $\varepsilon_i = \prod_{k=2}^i \frac{a_{k-1,k}}{a_{k,k-1}}$. In particular, $\varepsilon_i \, a_{i,i-1} = \varepsilon_{i-1} \, a_{i-1,i}$. 
In this proof, we omit the time dependence of the variables to simplify the notation, i.e. $S_t = S$ and $X_i(t) = X_i$. Then, using \eqref{model} and the definition of the $\varepsilon_i$s, we get that
\begin{align*}
\dot{S} & =  a_{10}^j X_1 -a_{11} X_1^2 -\sum_{i=2}^n \varepsilon_i a_{ii} X_i^2 - \sum_{i=2}^n \varepsilon_i a_{i0}^j X_i \\ 
& - a_{12} X_1 X_2 + \varepsilon_2 a_{21} X_1 X_2 + \sum_ {i=3}^n \varepsilon_i a_{i,i-1} X_{i-1} X_i -\varepsilon_{i-1} a_{i-1,i} X_{i-1} X_i \\
& \leq a_{10}^j X_1 -a_{11} X_1^2 -  \sum_{i=2}^n \varepsilon_i a_{i0}^j X_i.
\end{align*}
Let $\gamma = \max_{j \in E} a_{10}^j$ and $\varepsilon = \min_{j\in E} \min_{i=2,\ldots,n} \varepsilon_i a_{i0}^j$, then 
\begin{align*}
\dot{S} & \leq \gamma X_1 -a_{11} X_1^2 -  \varepsilon \sum_{i=2}^n X_i \\
& = (\gamma+\varepsilon) X_1 -a_{11} X_1^2 -  \varepsilon S \\
& \leq R - \varepsilon S
\end{align*}
for some constant $R>0$. 
If $S \geq \frac{R}{\varepsilon}$, then $\dot{S}\leq 0$ and if $S(0) \leq \frac{R}{\varepsilon}$, then $S(t) \leq \frac{R}{\varepsilon}$ for all $t \geq 0$. This makes the compact set
$$B := \left\{ x_1 + \sum_{i=2}^n \varepsilon_i x_i  \leq \frac{R}{\varepsilon} + 1\right\}$$
positively invariant under $\varphi^j$ (the flow induced by \eqref{functions_G^j}) for all $j \in E$.  Then $B \times E$ is positively invariant for \eqref{model} and for every initial condition $(x,j) \in \RR_+^n \times E$, with probability $1$ there exists a time $t\geq 0$ such that $Z^{(x,j)}(t) \in B\times E$.


\section{Sensitivity}  \label{Section_Sensibility}

In this part, we will further analyse the sensitivity of model \eqref{model}. First, we will see that if  all the vector fields are in some sense good for species $k$, then this species survives in the switched model. Moreover, if all environments are good for at most the $k$ first species, then species $k+1,\ldots,n$ become extinct in the PDMP setting. Note that, this property is not true in general, see \cite{BL16} for example.

In a second step, we will see that if we perturb  the coefficients $a_{ik}$, $k\neq 0$, in each environment, so that in environment $j$, $a_{ik}$ becomes $a_{ik}^j$, then for a sufficiently small perturbation, the noisy model still satisfies assumptions of Theorem \ref{theorem_stoch_pers} or of Theorem~\ref{theorem_extinction} if the non-disturbed one satisfies them. Note that in the following, we will see that it is not possible to calculate explicitly $\mathcal{I}^-_{k+1}$ for the perturbed model so that it is not possible to say something about the critical case (i.e. Theorem \ref{theorem_critical_case}) for the perturbed model. 

To do so, we first need to introduce some notations. Consider the more general situation where the vector fields $F^j$ are defined by \eqref{fucntions_F_i^j} but each coefficient depends on the environment, that is $a_{ik}$, $k\neq 0$, becomes $a_{ik}^j$.
For $a,b, m \in \NN$, $a \leq b$, we let $A_a^b (m)$ denote the set of permutations of $ \{a,a+1, \ldots, b-1, b\}$ that are product of $m$ disjoint transpositions of the form  $(i\; i+1)$ and we let
$$A_a^b := \bigcup_{m=0}^\infty \, A_a^b (m).$$
Remark that $A_a^b(m) = \emptyset$ if $m > \frac{b-a+1}{2}$. Now, we define for $1\leq k \leq n$ and $j\in E$,
\begin{align}\label{constant_delta}
\delta^j(k) := a_{10}^j \prod_{i=2}^k a_{i,i-1}^j - \sum_{m=2}^k a_{m0}^j \prod_{l=m+1}^k a_{l,l-1}^j \sum_{\alpha \in A_1^{m-1}} \prod_{i=1}^{m-1} a_{i,\alpha(i)}^j.
\end{align}
\begin{example}
\begin{align*}
A_1^3 &= \left\{ Id_{\{1,2,3\}}, (1 \; 2), (2 \;3) \right\},  \\
A_1^4 &= \left\{ Id_{\{1,2,3,4\}}, (1 \; 2), (2\; 3), (3 \; 4), (1 \; 2)(3 \; 4) \right\},\\
\delta^j(2) &= a_{10}^j\, a_{21}^j - a_{20}^j\, a_{11}^j,\\
\delta^j(3) &= a_{10}^j \,a_{32}^j \,a_{21}^j - a_{20}^j \,a_{32}^j \,a_{11}^j - a_{30}^j \,a_{22}^j \,a_{11}^j - a_{30}^j \,a_{12}^j \,a_{21}^j.
\end{align*}
\end{example}

We already know that $F^j_{|k}$ has a unique equilibrium $(\bar x_1^j,\ldots,\bar x_k^j)$ and moreover the proof of Proposition 2.4 of \cite{LV_degenerate} gives us an explicit formula for $\bar x_k^j$ ,
\begin{align}\label{Formule_explicite_zero}
\bar x_k^j \sum_{\alpha \in A_1^k} \prod_{i=1}^k a_{i,\alpha(i)}^j = \delta^j(k).
\end{align}
Since $F^\nu$ is still a Lotka-Volterra food chain, we define $\delta^\nu(k)$ in the same way with $a_{im}^j$ replaced by $a_{im}^\nu$ and the formula \eqref{Formule_explicite_zero} holds for $q^{*k}$, the equilibrium of $F^\nu_{|k}$. Another point is that by \eqref{Formule_explicite_zero} and the proof of Proposition~2.4 of \cite{LV_degenerate}, $F^j_{|k}$ has a unique positive equilibrium if and only if $\delta^j(k)>0$ (and the same holds for $F^\nu_{|k}$ and $\delta^\nu(k)$).

Let's go back to the case where the vector fields $F^j$ are defined by \eqref{fucntions_F_i^j}. Let
\begin{align}\label{Constant_Delta}
\Delta_k := \sum_{\alpha \in A_1^k} \prod_{i=1}^k a_{i,\alpha(i)}^j = \sum_{\alpha \in A_1^k} \prod_{i=1}^k a_{i,\alpha(i)}.
\end{align} 
The second equality follows from the fact that $a_{im}^j = a_{im}$ for $m\neq 0$.
Then, for $\mu ~\!\!\!\in~\!\!\! \PP_{erg}(\partial B \times E)$ such that $\mu(B_+^k \times E) = 1$, by the proof of Proposition \ref{delta_pos_iif_existe_p_i_st_sum_pi_mu_lmada_i_pos} we have that
\begin{align*}
\lambda_{k+1}(\mu) & = -a_{k+1,0}^\nu + a_{k+1,k}\,q^{*k}_k \\
& = \frac{1}{\Delta_k}\delta^\nu(k+1) \\
& = \frac{1}{\Delta_k}\left[ a_{10}^\nu \prod_{i=2}^{k+1} a_{i,i-1} - \sum_{m=2}^{k+1} a_{m0}^\nu \prod_{l=m+1}^{k+1} a_{l,l-1} \sum_{\alpha \in A_1^{m-1}} \prod_{i=1}^{m-1} a_{i,\alpha(i)} \right]\\
& =\frac{1}{\Delta_k} \sum_{j=1}^N \nu_j\left[ a_{10}^j \prod_{i=2}^{k+1} a_{i,i-1} - \sum_{m=2}^{k+1} a_{m0}^j \prod_{l=m+1}^{k+1} a_{l,l-1} \sum_{\alpha \in A_1^{m-1}} \prod_{i=1}^{m-1} a_{i,\alpha(i)} \right]\\
& = \frac{1}{\Delta_k} \sum_{j=1}^N \nu_j \, \delta^j(k+1) .
\end{align*}
The second equality follows from Lemma 5.1 of \cite{LV_degenerate}. Recall that by the proofs of the main theorems of this paper, the survival of species $k+1$ only depends on the sign of $\lambda_{k+1}(\mu)$. This works in the same manner for one environment $j\in E$ but with the sign of $\delta^j(k+1)$. Indeed, if $\delta^j(n)>0$, then by Proposition \ref{Accessibiliy_q*} environment $j$ causes all species to survive. Otherwise, by Lemma 5.1 of \cite{LV_degenerate} (or Proposition \ref{properties_equil}), if $\delta^j(i) >0$ then $\delta^j(i-1)>0$, this implies that there exists $m<n$ such that $\delta^j(m)>0$ and $\delta^j(m+1) \leq 0$. One can show that $(\bar x_1^j,\ldots, \bar x_m^j,0,\ldots,0)$ is globally asymptotically stable for the dynamic induced by $F^j$ where $(\bar x_1^j,\ldots, \bar x_m^j)$ is the positive equilibrium of $F^j_{|m}$ (this works in the same way as in \cite[Theorem 5.3.1]{HS98}). In this environment, species $1,\ldots,m$ survive while the others become extinct.

In other words, environment $j$ is favourable to species $k$ if and only if $\delta^j(k)>0$.
By the latter equation, this means that if all environments are good for species $k$ then the PDMP is also good for species $k$, meaning that with probability 1, species $k$ survives. Moreover, if all environments are good for at most $k$ species then the PDMP is bad for species $k+1,\ldots,n$.

Now, we will investigate the sensitivity of the model. Consider another PDMP Lotka-Volterra food chain defined by \eqref{model} with the same jump rates $b_{ij}$ and with also $N$ environments $\tilde F^j$ defined by \eqref{fucntions_F_i^j} but whose coefficients depend on the environment, that is 
\begin{align} \label{fucntions_F_i^j_v2}
\tilde F_i^j(x) := \left\{
     \begin{array}{ll}
      \tilde a_{10}^j - \tilde a_{11}^j x_1 - \tilde a_{12}^j x_2 								& i=1,\\
      - \tilde a_{i0}^j + \tilde a_{i,i-1}^j x_{i-1} - \tilde a_{ii}^j \, x_i - \tilde a_{i,i+1}^j x_{i+1} 	& i= 2, \ldots,n-1, \\
      - \tilde a_{n0}^j + \tilde a_{n,n-1}^j  x_{n-1} - \tilde a_{nn}^j \,x_n 					& i=n.
     \end{array}
     \right. 
\end{align}
Then there exist $\varepsilon>0$ such that 
\begin{align}\label{hypothese_robustesse}
|a_{im}^j - \tilde a_{im}^j| \leq \varepsilon \qquad \forall i,m \text{ and } \forall j \in E.
\end{align}

We call the {\it basic model} for the one with the environments $F^j$ and the {\it noisy model} for the one with the environments $\tilde F^j$.

%
\begin{proposition}\label{lambda_same_sign}
Let $\mu$ be an ergodic probability measure for the basic model and let $\mu_\varepsilon$ be another ergodic probability measure for the noisy one. Suppose that $I_\mu = I_{\mu_\varepsilon} = \{1 , \ldots, k\}$ so that $\mu(B_+^k \times E) = \mu_\varepsilon(B_+^k \times E) = 1$, then for $\varepsilon$ sufficiently small, $\lambda_{k+1} (\mu)$ and $\lambda_{k+1}(\mu_\varepsilon)$ have the same sign.
\end{proposition}

\prf 
By Proposition \ref{Invasion_rate_0}, we have that
$$\lambda_i(\mu_\varepsilon) =0  \qquad i=1,\ldots,k,$$
that is
$$\left\{
     \begin{array}{ll}
     0 = \tilde a_{10}^\nu - \sum_{j=1}^N \tilde a_{11}^j \int x_1 d\mu_\varepsilon^j - \sum_{j=1}^N \tilde a_{12}^j \int x_2 d\mu_\varepsilon^j,		&\\
     0 = - \tilde a_{i0}^\nu + \sum_{j=1}^N \tilde a_{i,i-1}^j \int x_{i-1} d\mu_\varepsilon^j - \sum_{j=1}^N \tilde a_{ii}^j \int x_i d\mu_\varepsilon^j - \sum_{j=1}^N \tilde a_{i,i+1}^j \int x_{i+1} d\mu_\varepsilon^j 	& i\neq
      1,k, \\
     0 = - \tilde a_{k0}^\nu + \sum_{j=1}^N \tilde a_{k,k-1}^j \int x_{k-1} d\mu_\varepsilon^j -\sum_{j=1}^N \tilde a_{kk}^j \int x_k d\mu_\varepsilon^j.	&
     \end{array}
\right. $$
Let $\tilde x_i^\varepsilon := \sum_{j=1}^N \int x_i d\mu_\varepsilon^j$, then by \eqref{hypothese_robustesse}
$$\left\{
     \begin{array}{ll}
     a_{10}^\nu + \varepsilon \geq  ( a_{11} -\varepsilon) \, \tilde x_1^\varepsilon + ( a_{12} -\varepsilon) \, \tilde x_2^\varepsilon,		&\\
     0 \leq - (a_{i0}^\nu - \varepsilon) +  ( a_{i,i-1} +\varepsilon) \, \tilde x_{i-1}^\varepsilon - ( a_{ii} - \varepsilon) \, \tilde x_i^\varepsilon - ( a_{i,i+1} - \varepsilon) \, \tilde x_{i+1}^\varepsilon 	& i\neq
      1,k, \\
     0 \leq - (a_{k0}^\nu - \varepsilon) + ( a_{k,k-1} +\varepsilon) \, \tilde x_{k-1}^\varepsilon -( a_{kk} -\varepsilon) \, \tilde x_k^\varepsilon.	&
     \end{array}
\right. $$
Then it is straightforward that
\begin{align}\label{bound_up_x_k_epsilon}
\tilde x_k^\varepsilon \leq f_k(\varepsilon) = \frac{\delta^\varepsilon(k)}{\Delta_k^\varepsilon}
\end{align}
where $\delta^\varepsilon(k)$ and $\Delta_k^\varepsilon$ are defined by \eqref{constant_delta} and \eqref{Constant_Delta} respectively but with $a^j_{10}$ replaced by $a^\nu_{10}+ \varepsilon$, $a^j_{i0}$ by $a^\nu_{i0}- \varepsilon$ $i\neq 1$, $a^j_{i,i-1}$ by $a_{i,i-1}+ \varepsilon$ and $a^j_{ii}, a^j_{i,i+1}$ by $a^j_{ii} - \varepsilon, a^j_{i,i+1}-\varepsilon$.

In the same manner we obtain that
\begin{align}\label{bound_down_x_k_epsilon}
\tilde x_k^\varepsilon \geq f_k(-\varepsilon) = \frac{\delta^{-\varepsilon}(k)}{\Delta_k^{-\varepsilon}}
\end{align}
where $\delta^{-\varepsilon}(k)$ and $\Delta_k^{-\varepsilon}$ are defined as before with $\varepsilon$ replaced by $-\varepsilon$. Note that $ \lim_{\varepsilon \to 0} f_k(\varepsilon) =\lim_{\varepsilon \to 0} f_k(-\varepsilon) = q^{*k}_k$.

Then if we let $g_{k+1}(\varepsilon) :=-(a_{k+1,0}^\nu - \varepsilon) + (a_{k+1,k} +\varepsilon) f_k(\varepsilon)$,
\begin{align*}
g_{k+1}(-\varepsilon) \leq \lambda_{k+1}(\mu_\varepsilon) \leq g_{k+1}(\varepsilon).
\end{align*}
In particular, $\lim_{\varepsilon \to 0} g_{k+1}(\varepsilon) = \lambda_{k+1}(\mu)$ which concludes the proof.

{\hfill$\square$\bigbreak}

As a result, if the noisy model satisfies Assumption \ref{Assumption_bio} and the basic one satisfies assumptions of Theorem \ref{theorem_stoch_pers} or Theorem \ref{theorem_extinction} then the conclusions of these theorems are also true for the noisy model.

Indeed, suppose first that the basic model satisfies assumptions of Theorem \ref{theorem_stoch_pers}. Then for $\varepsilon>0$ sufficiently small, by the previous proposition and Proposition \ref{delta_pos_iif_existe_p_i_st_sum_pi_mu_lmada_i_pos}, condition \eqref{Criter_alpha_i_lambda_i} of Theorem \ref{Criter_thm_stoc_persi} is satisfied. Consider $\tilde G^\nu$ the average vector field of the noisy model, then
$$a_{ik}^\nu - \varepsilon \leq \tilde a_{ik}^\nu \leq a_{ik}^\nu + \varepsilon \qquad \forall \, i,k.$$
Like in the proof of Proposition \ref{lambda_same_sign}, for $\varepsilon>0$ sufficiently small, as $F^\nu$ has a positive equilibrium, so does $\tilde F^\nu$. Then by Proposition \ref{Accessibiliy_q*}, Assumption \ref{Assumption_bio}, and Proposition \ref{Hormander_at_q*}, this equilibrium is accessible and verifies the strong Hörmander condition. By Theorem \ref{Criter_thm_stoc_persi}, the same conclusions as for basic model apply for the noisy one.

Now suppose that the basic model satisfies assumptions of Theorem \ref{theorem_extinction}. By Proposition \ref{lambda_same_sign}, for sufficiently small $\varepsilon>0$, the proof of Theorem \ref{theorem_extinction} works in the same way for the noisy model. The only thing that needs to be clarified is the rate of extinction of species $k+1$, $\mathcal{I}_{k+1}^{\varepsilon-}$. Indeed, we cannot have an explicit formula for this since we are not able to calculate $\lambda_{k+1}(\mu_\varepsilon)$ explicitly for $\mu_\varepsilon$ an ergodic probability measure for the noisy model such that $I_{\mu_\varepsilon} = \{1,\ldots,k\}$.

However, with the help of the bounds of $\tilde x_{k+1}^\varepsilon$ given by \eqref{bound_up_x_k_epsilon} and \eqref{bound_down_x_k_epsilon}, we can have an estimate of it
$$\frac{\delta^{-\varepsilon}(k+1)}{\Delta_k^{-\varepsilon}} \leq \mathcal{I}_{k+1}^{\varepsilon-} \leq \frac{\delta^\varepsilon(k+1)}{\Delta_k^\varepsilon}.$$
So that the extinction part works in the same way for the noisy model as for the basic one.

Note that it also works the other way round. Indeed, let us consider a PDMP Lotka-Volterra food chain with vector fields $\tilde F^j$ defined by \eqref{fucntions_F_i^j_v2} which satisfies Assumption \ref{Assumption_bio}. In order to have the conclusions of Theorem \ref{theorem_stoch_pers} or of Theorem \ref{theorem_extinction}, it suffices to find $\varepsilon>0$ sufficiently small and coefficients $a_{ik}$, $k\neq 0$ such that
$$|a_{ik}-\tilde a_{ik}^j | < \varepsilon \qquad \forall i,k \text{ and } j\in E$$
and such that the PDMP with vector fields $F^j$ with coefficients $a_{ik}$, $k\neq 0$ and $a_{i0}^j$ satisfies assumptions of Theorem \ref{theorem_stoch_pers} or Theorem \ref{theorem_extinction} and Proposition \ref{lambda_same_sign}.


\section*{Acknowledgments}
This work was supported by the SNF grants 200020-196999. I thank Michel Benaïm and Edouard Strickler for useful discussions and reviews, and Oliver Tough for his kind comments on the introduction.

I thank the two anonymous referees for their useful comments on the presentation of the main results.  In particular, it helped me to realise that the critical case was not being addressed, which is now the case.


\bibliographystyle{amsplain}
\bibliographystyle{imsart-nameyear}
\bibliographystyle{nonumber}
\bibliography{PDMP_LV_food_chain_Biblio}
\end{document}